\documentclass[11pt,a4paper]{article}
\usepackage[francais,english]{babel}
\usepackage[latin1]{inputenc}
\usepackage{hyperref}
\usepackage{latexsym,amssymb,amsmath,amsfonts,delarray,epsfig,eucal,graphicx,hyperref}
\usepackage{amsmath}
\usepackage{amssymb}
\usepackage[T1]{fontenc}
\usepackage{multicol}
\usepackage{calc}
\usepackage[all]{xy}
\usepackage{amsthm}

\theoremstyle{plain}
\newtheorem{theoreme}{\sc{Theorem}}[section]
\newtheorem{thm}{\sc{Theorem}}
\newtheorem{proposition}[theoreme]{Proposition}
\newtheorem{corollaire}[theoreme]{Corollary}
\newtheorem{lemme}[theoreme]{Lemma}

\theoremstyle{definition}
\newtheorem{definition}[theoreme]{Definition}

\theoremstyle{remark}

\newtheorem{remarque}[theoreme]{Remark}

\title{Relative Zariski Open Objects}
\author{Florian Marty\\fmarty@math.ups-tlse.fr}
\date{Université Toulouse III - Laboratoire Emile Picard}
\begin{document}
\maketitle
\setcounter{secnumdepth}{-1}
\section{Abstract}
\footnotesize{In \cite{8}, Bertrand To\"en and Michel Vaqui\'e define a scheme theory for a closed monoidal category $(\mathcal{C},\otimes,1)$. One of the key ingredients of this theory is the definition of a Zariski topology on the category of commutative monoids in $\mathcal{C}$. The purpose of this article is to prove that under some hypotheses, Zariski open subobjects of affine schemes can be classified almost as in the usual case of rings $(Z-mod,\otimes,Z)$. The main result states that for any commutative monoid $A$, the locale of Zariski open subobjects of the affine scheme $Spec(A)$ is associated to a topological space whose points are prime ideals of $A$ and open subsets are defined by the same formula as in rings. As a consequence, we compare the notions of scheme over $\mathbb{F}_{1}$ of \cite{10} and \cite{8}.} 
\setcounter{secnumdepth}{2}
\small
\renewcommand{\theenumi}{\roman{enumi}}
\tableofcontents
\setcounter{secnumdepth}{-1}
\section{Introduction}
In \cite{8}, Bertrand Toën and Michel Vaquie define a scheme theory for a closed monoidal category $(\mathcal{C},\otimes,1)$. This theory recovers the classical notion of scheme when $(\mathcal{C},\otimes,1)=(\mathbb{Z}-mod,\otimes,\mathbb{Z})$. One of its key ingredients is the definition of a Zariski topology on the category of commutative monoids of $\mathcal{C}$. This topology is used to glue affine schemes, which by definition are commutative monoids in $\mathcal{C}$, in order to obtain the general notion of scheme over $\mathcal{C}$. The aim of this article is to study in more detail this Zariski topology.\vskip 4pt
The definition of Zariski open object of \cite{8} differs from the usual definition that could not be generalized to a relative context. In fact, a theorem in SGA asserts that the usual notion of Zariski open is equivalent to the notion of \textit{finitely presented flat epimorphism} which are well defined in the relative context.\\
The usual notions of image, generators, ideal, prime ideal or localization of a monoid have generalisations to this relative setting. For instance, for a monoid $A$, the localization by an element $f$ of its underlying set $A_{0}:=Hom_{\mathcal{C}}(1,A)$, denoted by $A_{f}$, makes sense. Furthermore, they are the first examples of relative Zariski open objects and they play a fundamental role in the main theorems:
\begin{thm}\ref{basis}
Let $A$ be a commutative monoid. Zariski open objects $(A_{f})_{f\in A_{0}:=Hom_{\mathcal{C}}(1,A)}$ form a basis of open objects of the Zariski topology of $Spec(A)$.
\end{thm}
\begin{thm}\ref{spac}
Let $A$ be a commutative monoid. The poset of Zariski open subobjects of $X:=Spec(A)$ is the locale associated to the topological space $Ouv(X)$ whose points are prime ideals of $A$ and closed sets are given, for an ideal $q$ of $A$, by $V(q):=\{p\;st\;q\subset\;p\}$.
\end{thm}
These theorems are well known for $\mathcal{C}=\mathbb{Z}-mod$. In a relative context, we prove that there is a generalisation of the sober topological space associated to a ring. For a commutative monoid $A$ in $\mathcal{C}$ this construction is not natural. The theory of enriched sheaves of Borceux and Quinteiro plays a fundamental rôle here. In particular, the notion of enriched Grothendieck topology will provide us with a notion of filter of ideals. Such filters will be called Gabriel filters. We will prove that a Zariski open subobject has an associated Gabriel filter which is more regular in the sense that it is associated to one ideal. These Gabriel filters are called locally primitive Gabriel filters. For a commutative monoid $A$, we prove that there is a contravarient equivalence of poset between the poset of locally primitive Gabriel filters of $A$ and the poset of Zariski open subobjects of $Spec(A)$. We conclude by proving that these Gabriel filters are characterized by their subsets of prime ideals.\\
Finally, we use these results in the relative context $\mathcal{C}=(Ens,\times,\mathbb{F}_{1})$ to compare the notions of $\mathbb{F}_{1}-scheme$ of \cite{8} and \cite{10} and to prove that the topological spaces associated are homeomorphic.\vskip 8pt

\noindent \textbf{Acknowledgements:} I wish to thank Bertrand Toen for his helpfull comments and suggestions on this work. I also thank Joseph Tapia for the clearness of his comments on category theory and Michel Vaquie for a usefull feedback on this work, in particular on a (necesarilly vicious) finitude problem.
 
\section{Preliminaries}
All along this work $(\mathcal{C},\otimes,1)$ is a closed symmetric monoidal category which is complete, cocomplete, locally finitely presentable, regular in the sense of Barr and has a compact(ie finitely presented) unit. Recall that $\mathcal{C}$ is regular in the sense of Barr if the pullback of a regular epimorphism (i.e. a cokernel) is also a regular epimorphism and $\mathcal{C}$ is locally finitely presentable if the Yoneda functor $i:\mathcal{C}\rightarrow pr(\mathcal{C}_{0})$ is fully-faithfull, where $\mathcal{C}_{0}$ is the full subcategory of finitely presented objects. We assume futhermore that $\mathcal{C}_{0}$ is stable under tensor products and that $Hom_{\mathcal{C}}(1,-)$ preserves regular epimorphisms (we say that $\mathcal{C}$ respects images). A category $\mathcal{C}$ verifying these hypothese will be called a relative context.\vskip 4pt
The functor $i$ has a left adjoint $K$ defined, for $F\in Pr(\mathcal{C}_{0})$ by $K(F):=Colim_{(k,f)\in \mathcal{C}_{0}^{F}}k$ where $\mathcal{C}_{0}^{F}$ is the category whose objects are couples $(k,f)$, $k\in \mathcal{C}_{0}$, $f\in F(k)$ and morphisms from $(k,f)$ to $(k',f')$ are morphisms $h:k\rightarrow k'$ such that $f'\circ h=f$.\vskip 4pt 
The hypotheses have three first important consequences. The functor $i$ reflects epimorphisms, filtered colimits are exact in $\mathcal{C}$ and the collection of subobjects of a given object $X$, ie of isomorphism classes of couple $(Y,j)$ where $\xymatrix{Y\ar[r]^{j}&X}$, $j$ is a monomorphism, is a set.
\vskip 4pt
The functors $Hom_{\mathcal{C}}(1,-)$, denoted $(-)_{0}:X\rightarrow X_{0}$ and $Hom_{\mathcal{C}}(k,-)$, denoted $(-)_{k}:X\rightarrow X_{k}$, $k$ in $\mathcal{C}_{0}$ from  $\mathcal{C}$ to Set will be called respectively \textit{underlying set functor} and \textit{weak underlying set functor}. For $X\in \mathcal{C}$, the object $X_{0}$ is the underlaying set of $X$. If there is no confusion, we may forget the word $weak$ for the other functors. An element of $X$ is an element of one of its weak underlying sets $X_{k}$ or of $X_{0}$. If $f$ is such an element, we will write $f\in X$.\vskip 4pt A relative context is called \textit{strong} when the functor $(-)_{0}$ reflects isomorphisms i.e. when the unit $1$ is generating.\vskip 4pt
For $X,Y\in \mathcal{C}$, let $r_{X}$, $l_{X}$, $s_{(X,Y)}$ and $a_{X,Y,Z}$ denote respectively the natural isomorphisms $X\otimes 1\rightarrow X$, $1\otimes X\rightarrow X$, $X\otimes Y\rightarrow Y\otimes X$, and $(X\otimes Y)\otimes Z\rightarrow X\otimes(Y\otimes Y)$. The Mac Lane coherence theorem (cf \cite{7}) asserts that any monoidal category is equivalent to one in which morphisms $r,l$ and $a$ are identities, so they can be considered as identities in all demonstrations concerning properties of $\mathcal{C}$ that are stable under monoidal equivalences.\vskip 4pt
In such a category, there exists a notion of commutative monoid and for a given commutative monoid $A$, of $A$-module (see \cite{7} or \cite{8}).  For a commutative monoid $A$ and an $A$-module $M$, let $m_{A}$, $i_{A}$, $\mu_{M}$ denote respectively the multiplication $A\otimes A\rightarrow A$, the unit morphism $1\rightarrow A$  and the multiplication $M\otimes A\rightarrow M$. It is known that the categories of modules over a monoid $B$ are also symmetric monoidale categories. The notations for structural morphisms will then be used as well for the category $\mathcal{C}$ as for any category $B-mod$. $Comm(\mathcal{C})$ will refer to the category of commutative monoids in $\mathcal{C}$ and for $A\in Comm(C)$, $A-mod$ will refer to the category of $A$-modules.\vskip 4pt
In fact, all the hypotheses on our basis category will be true in any of its sub-categories of modules. Let $A$ be a commutative monoid in $\mathcal{C}$. The category $A-mod$ is a relative context. Let us give some details on this fact. It's well known that $A-mod$ is closed, complete and cocomplete. Moreover, as the forgetfull functor $A-mod\rightarrow \mathcal{C}$ commutes with small limits and colimits, the regularity in the sense of Barr, the compactness of the unit $A$ and the fact that $A-mod$ respects images are easy to prove in $A-mod$. The most difficult is to prove that the full sub-category of finitely presentable objects of $A-mod$, denoted $A-mod_{0}$ is a generating subcategory. It can be proved for example using Ind-objects.\vskip 4pt
In a relative context, the category of commutative monoids has also some interresting properties. In particular, it can be proved that it is locally finitely presentable. As the categories of modules are also relative contexts, these properties are true in any category of algebras in $\mathcal{C}$. Recall that, for the monoid $A$, the categories of algebras over $A$, denoted $A-alg$ is $Comm(A-mod)$. It is also equivalent to $A/Comm(\mathcal{C})$. The category of finitely presentable objects in $A-mod$ will be denoted $A-mod_{0}$.
\vskip 4pt
Let us now recall some well known properties of monoidal categories. A pushout in $A-alg$ is a tensor product in the sense that for commutative monoids $B,C\in A-alg$,  $B\otimes_{A}C\backsimeq B\coprod_{A}C$.
There are two adjunctions:
\begin{center}
$\xymatrix{\mathcal{C}\ar@<1ex>[r]^-{(-\otimes A)}&A-mod\ar@<1ex>[l]^-{i}}$\;\;\;\;\;\;\;\;\;\;$\xymatrix{\mathcal{C}\ar@<1ex>[r]^-{L}&Comm(\mathcal{C})\ar@<1ex>[l]^-{i}}$
\end{center}
where the forgetfull functor $i$ is a right adjoint and the \textit{free associated monoid} functor $L$ is defined by $L(X):=\coprod_{n\in \mathcal{N}}X^{\otimes n}/S_{n}$ ($S_{n}$ refers to the symmetric group). In these adjunctions, $\mathcal{C}$ can be replaced by $B-mod$ for $B\in Comm(\mathcal{C})$. Let $\varphi$ (resp $\varphi_{B}$) and $\psi$ (resp $\psi_{B}$) denote these adjunctions for the category $\mathcal{C}$ (resp $B-mod$). For $X\in \mathcal{C}$ and $M\in A-mod$, $\varphi:Hom_{\mathcal{C}}(X,M)\rightarrow Hom_{A-mod}(X\otimes A,M)$ is easy to describe :
\begin{center}
$\varphi: f\rightarrow \mu_{M}\circ Id_{A}\otimes f$\vskip 3pt$\varphi^{-1}:g\rightarrow g\circ (Id_{X}\otimes i_{A})\circ r_{X}^{-1}$
\end{center}
\setcounter{secnumdepth}{2}

\section{Algebraic Structures}
\subsection{Images and Generators}
We give and describe now two fundamental notions in a relative context $\mathcal{C}$. These notions will be particularly important in the relative contexts $A-mod$ for a commutative monoid $A$.  
\subsubsection{Images}
\begin{definition}Let $u:X\rightarrow Y$ be a morphism in $\mathcal{C}$.
\begin{list}{$\triangleright$}{}
\item Let $\mathcal{C}_{Y}$ denote the full sub-category of $\mathcal{C}/Y$ consisting of objects $Z$ such that the morphism $Z\rightarrow Y$ is a monomorphism. A subobject of $Y$ is an isomorphism class of objects in $\mathcal{C}_{Y}$. The category of subobjects of $Y$ will be denoted $sub(Y)$. 
\item  A subobject $Z$ of $Y$ \textit{contains the image of $u$} if $u$ factors throught the inclusion $Z\hookrightarrow Y$. We denote by $sub_{u}Y$ the category whose objects are subobjects of $Y$ \textit{containing the image of $u$} and whose morphisms are monomorphisms.
\item The image of $u$, denoted $Im(u)$ is
\begin{center}
$Im(u):=Lim_{sub_{u}Y}Z$.
\end{center}
\item Let $|u|$ denote the object $|u|:=Coker(\xymatrix{X\times_{Y}X\ar@<1ex>[r]\ar@<-1ex>[r]&X})$.
\end{list}
\end{definition}
The existence of the limit is due to the smallness of the category $sub_{u}Y$. 
\begin{lemme}Let $u:X\rightarrow Y$ be a morphism in $\mathcal{C}$. 
\begin{list}{$\diamond$}{}
\item If $u$ is a monomorphism, there is an isomorphism $X\backsimeq Im(u)$. 
\item Let $v:Y\rightarrow Z$ be in $\mathcal{C}$. There is an isomorphism $Im(v\circ u)\backsimeq Im(v\circ u')$ where $u'$ is the monomorphism $u':Im(u)\rightarrow Y$.
\item There is an isomorphism $Im(u)\backsimeq |u|$.
\end{list}
\end{lemme}
Proof
\begin{list}{-}{}
\item The isomorphism class associated to the object $X$ is initial in $sub_{u}Y$.
\item This isomorphism comes from the equality $v\circ u=v\circ u'\circ u''$ where $u''$ is the morphism $u'':X\rightarrow Im(u)$. And by the universal property of the image, $Im(v\circ u'\circ u'')=Im(v\circ u')$. 
\item There is a unique natural morphism from $|u|$ to any subobject of $Y$ \textit{containing the image of $u$} and thus from $|u|$ to $Im(u)$. We need then to prove that $|u|$ is a subobject of $Y$ to achieve the proof. Let us prove first that $X\times_{Y}|u|\backsimeq X$. As a relative context is regular in the sense of Barr, 
\begin{center}
$|u|\times_{Y}X\backsimeq Coker(\xymatrix{X\times_{Y}X\times _{Y}X\ar@<1ex>[r]\ar@<-1ex>[r]&X\times_{Y}X})\backsimeq |p|$
\end{center}
where the natural morphism $p:X\times_{Y}X\rightarrow X$ is split by a morphism $s$, i.e. $p\circ s = Id$. Let $r,q$ and $t$ denote respectively the morphism $X\times_{Y}X\rightarrow |p|$, the morphism $|p|\rightarrow X$ and the morphism $r\circ s$. More precisely, there is a diagram
\begin{center}
$\xymatrix{X\times_{Y}X\ar[rr]^{r}\ar@<1ex>[rd]^{p}&&|p|\ar@<1ex>[ld]^{q}\\&X\ar@<1ex>[lu]^{s}\ar@<1ex>[ru]^{t}}$
\end{center}
in which $q\circ t= q\circ r\circ s=p\circ s= Id$. Reciprocally, let us prove that $t\circ q= Id$. As $r$ is an epimorphism, it is enought to prove $t\circ q\circ r = r$. As $t\circ q\circ r = r\circ s\circ q\circ r$ and $q\circ r=p$ (by construction of $q$), then $t\circ q\circ r=r\circ s\circ p=r$.\vskip 2pt
Finally
\begin{center}
$|u|\times_{Y}|u|\backsimeq Coker(\xymatrix{X\times_{Y}X\times _{Y}|u|\ar@<1ex>[r]\ar@<-1ex>[r]&X\times_{Y}|u|})\backsimeq |u|$
\end{center}
thus $|u|$ represents a subobject of $Y$.   
\end{list} 
\vskip 2pt
$\hfill\blacklozenge$
\begin{corollaire}
Let $u:X\rightarrow Y$ be a morphism in $\mathcal{C}$. The morphism $u_{0}:X_{0}\rightarrow Y_{0}$ is surjective.
\end{corollaire}
Proof\\As $\mathcal{C}$ respects images, the functor $(-)_{0}$ commutes with regular epimorphim and thus send them to epimorphisms in the category of sets, i.e. to surjective morphisms.
\vskip 2pt
$\hfill\blacklozenge$
\begin{lemme}
Let $u:A\rightarrow B$ be a morphism in $Comm(\mathcal{C})$. Its image in $\mathcal{C}$ is isomorphic (in $\mathcal{C}$) to its image in $Comm(\mathcal{C})$.
\end{lemme}
Proof\\
Let $Im(u)$ be its image in $\mathcal{C}$. We provide it with a structure of commutative monoid. Let $u'$ denote the morphism $A\rightarrow Im(u)$. The unit morphism of $Im(u)$ is $u'\circ i_{A}$. Let us now define a multiplication morphism. Recall that the forgetfull functor from $A-mod$ to $\mathcal{C}$ commutes with small colimits. Thus $Im(u)$ has a natural structure of $A$-module. As $\mathcal{C}$ is closed, functors $X\otimes -$, $X\in \mathcal{C}$ commute with small colimits and preserve epimorphisms. Thus
\begin{center}
$Im(u)\otimes Im(u)\backsimeq Coker(\xymatrix{A\times_{B}A\otimes Im(u)\ar@<1ex>[r]\ar@<-1ex>[r]&A\otimes Im(u)})$\vskip 2pt $Im(u)\otimes A\backsimeq Coker(\xymatrix{A\times_{B}A\otimes A\ar@<1ex>[r]\ar@<-1ex>[r]&A\otimes A})$
\end{center}
Moreover, the morphism $Id\otimes u':A\times_{B}A\otimes A\rightarrow A\times_{B}A\otimes Im(u)$ is an epimorphism. There is a commutative diagram
\begin{center}
$\xymatrix{A\times_{B}A\otimes A\ar@<1ex>[r]\ar@<-1ex>[r]\ar[d]_{Id\otimes u'}&A\otimes A\ar[d]^{Id_{A}\otimes u'}\\A\times_{B}A\otimes Im(u)\ar@<1ex>[r]\ar@<-1ex>[r]&A\otimes Im(u)\ar[r]_{\mu_{Im(u)}}&Im(u)}$
\end{center} 
where $Id_{A}\otimes u'$ equalizes the two top morphisms. Thus $\mu_{Im(u)}$ equalizes the two morphisms from $A\times_{B}A\otimes A$ to $A\times Im(u)$ and as $Id\otimes u'$ is an epimorphism, $\mu_{Im(u)}$ equalizes the two morphisms of the coker $Im(u)\otimes Im(u)$. Finally $\mu_{Im(u)}$ has a lift $m_{Im(u)}:Im(u)\otimes Im(u)\rightarrow Im(u)$. The commutativity of the structure diagrams of this monoid is clear.
\vskip 2pt
$\hfill\blacklozenge$

\subsubsection{Generators}
\begin{definition}
Let $X$ be in $\mathcal{C}$ and $F=(f_{i})_{i\in I}$ be a family of elements of $X$. Let $sub_{F}X$ denote the full subcategory of $sub(X)$ whose objects are subobjects $Y$ of $X$ containing the family $F$. The object generated by $F$ is 
\begin{center}
$X_{F}:=Lim_{sub_{F}X}Y$
\end{center}
\end{definition}
\begin{definition}
Let $X$ be in $\mathcal{C}$ and $F=(f_{i})_{i\in I}$ be a family of elements of $X$. The presheaf associated to $F$, denoted $F'$, is defined by
\begin{center}
$F':k\rightarrow F(k):=\{f\in X_{k}\;st\;\exists\;h:k\rightarrow k'\;and\;g\in F\cap C_{k'}\;st\;f=g\circ h\}$
\end{center}
\end{definition}
It is a presheaf by construction. Moreover, the family $\coprod_{k}F'(k)$, still denoted $F'$, generates the same object, i.e. $X_{F'}\backsimeq X_{F}$.
\begin{lemme}
Let $X$ be in $\mathcal{C}$, $F=(f_{i})_{i\in I}$ be a family of elements of $X$ and $F'$ be its associated presheaf. Then the object $X_{F}$ is isomorphic in $\mathcal{C}$ to the image of the natural morphism
\begin{center}
$p:K(F'):=Colim_{\mathcal{C}_{0}^{F'}}k\rightarrow X$
\end{center}
where $K$ is the left adjoint of the Yoneda functor $i:\mathcal{C}\rightarrow Pr(\mathcal{C}_{0})$.
\end{lemme}
Proof\\Considering the adjunction between $K$ and $i$, the category $sub_{p}X$ is clearly isomorphic to $sub_{F}X$.  
\vskip 2pt
$\hfill\blacklozenge$
\vskip 2pt
Let us now give the properties of an object generated by a family.
\begin{lemme}\label{prgen}
Let $X$ be in $\mathcal{C}$ and $F$ be a family of elements of $X$.
\begin{enumerate}
\item If $Y,Y'$ are two objects of $sub(X)$ generated by $F$, they are isomorphic.
\item If $(-)_{0}$ reflect isomorphisms, $X_{0}$ is a generating family of $X$. 
\item Let $u:X\rightarrow Y$ be in $\mathcal{C}$. The family $u_{*}(F)$ generate $Im(u)$.
\item The family $F$ generates $X$ if and only if it is an epimorphic family of morphisms.
\item A morphism is an epimorphism if and only if it preserves generating families.
\item Let $\mathcal{G}$ be the diagram of finite subsets of $F$ whose morphisms are inclusions. There is an isomorphism $X_{F}\backsimeq Colim_{G\in \mathcal{G}}X_{G}$.
\item If $X$ is generated by $F$ and finitely presented in $\mathcal{C}$ ($Hom_{\mathcal{C}}(X,-)$ commute with filtered colimits), it is generated by a finite family included in $F$.
\item Let $A,B$ be in $\mathcal{C}$. The family $i(A)\times i(B)$ generates $A\otimes B$.  
\item Let $u:A\rightarrow B$ and $C$ be in $\mathcal{C}$. The image of $u\otimes Id_{C}$ is generated by $i(Im(u))\otimes i(C)$.
\end{enumerate}
\end{lemme}
Proof
\begin{enumerate}
\item Clear
\item Let us assume that $(-)_{0}$ reflect isomorphisms. Let $Y$ be the subobject of $X$ generated by the family of morphisms included in $X_{0}$. There is an isomorphism $X_{0}\backsimeq Y_{0}$, thus $X\backsimeq Y$.
\item The elements of $u_{*}(F)$ are clearly elements of $Im(u)$. Let $Im(u)_{F}$ be the object generated by $u_{*}(F)$. The natural morphism $F\rightarrow Im(u)_{F}$ lifts to a morphism $X\rightarrow Im(u)_{F}$ which factors $u$. Thus there is a monomorphism $Im(u)\rightarrow Im(u)_{F}$. As $i:\mathcal{C}\rightarrow Pr(\mathcal{C})$ is fully faithfull, there is an isomorphism $Im(u)\backsimeq Im(u)_{F}$.
\item A generating family is clearly epimorphic. Reciprocally, if $F$ is epimorphic, let $F'$ be its associated presheaf. There is an isomorphism for all $Y\in \mathcal{C}$, $Hom_{\mathcal{C}}(X,Y)\backsimeq Hom_{\mathcal{C}}(K(F'),Y)$. Thus, by Yoneda, $X\backsimeq K(F')\backsimeq X_{F}$.
\item Clear, epimorphisms preserves epimorphic families by composition.
\item We prove first that $Colim_{G\in \mathcal{G}}X_{G}$ is a subobject of $X$. The Yoneda functor $i:\mathcal{C}\rightarrow Pr(\mathcal{C}_{0})$ commutes with filtered colimits thus we compute this colimit levelwise in $Pr(\mathcal{C}_{0})$. More precisely, there are isomorphisms $colim_{G\in \mathcal{G}}(X_{G})_{k}\backsimeq Colim_{G\in \mathcal{G}}((X_{G})_{k})$. For all $k\in \mathcal{C}_{0}$, the sets $(X_{G})_{k}$ are subsets of $X_{k}$, therefore their filtered colimit is a subset of $X_{k}$. This implies that $Colim_{G\in \mathcal{G}}X_{G}$ is a subobject of $X$ in $\mathcal{C}$. Finally, any element of $i(X)$ belongs to a $i(X_{G})$ hence to their colimit and $X\backsimeq Colim_{G\in \mathcal{G}}X_{G}$.
\item  The object $X$ is finitely presentable and isomorphic to $Colim_{G\in \mathcal{G}}X_{G}$ thus there exists a finite set $G$ such that this isomorphism factors throught $X_{J}$. The morphism $X_{J}\rightarrow X$ is then a monomorphism and a split epimorphism hence an isomorphism.
\item We write $A\otimes B\backsimeq K(i(A))\otimes K(i(B))$, thus
\begin{center}
$A\otimes B\backsimeq colim_{(k,f),(k',f')\in \mathcal{C}^{i(A)}_{0}\times \mathcal{C}^{i(B)}_{0}}k\otimes k'\backsimeq K(F)$
\end{center}
where $F:=i(A)\times i(B)$.
\item The morphism $u\otimes Id_{C}$ factors through $Im(u)\otimes C\rightarrow B\otimes C$ and $Im(u)\otimes C$ is generated by $i(Im(u))\times i(C)$. The morphism $i(A)\times i(C)\rightarrow i(Im(u))\times i(C)$ is surjective thus the image of $u\otimes Id_{C}$ is isomorphic to the image of $Im(u)\otimes C\rightarrow B\otimes C$ which is generated by $i(Im(u))\times i(C)$. 
 
\end{enumerate}
\vskip 2pt
$\hfill\blacklozenge$

\subsection{Ideals and Multiplications}

\subsubsection{Ideals}
\begin{definition}
Let $A$ be a commutative monoid and $M$ be an $A$-module. A sub-$A$-module of $M$ is a subobject of $M$ in the category $A-mod$. An ideal of $A$ is a sub-$A$-module of $A$.
\end{definition}
In the rest of the paper, for any ideal $q$, an element of the isomorphism class will be chosen and also called $q$. The corresponding monomorphism will be denoted $j_{q}:q\rightarrow A$.
\begin{lemme}\label{r0}Let $A$ be in $Comm(\mathcal{C})$, $M$ be an $A$-module and $q$ be an ideal of $A$.
\begin{list}{$\diamond$}{}
\item The object $A_{0}$ is a monoïd in $(Set,\times,*)$, the object $M_{0}$ is a $A_{0}$-module in $(Set,\times,*)$.
\item The object $q_{0}$ is a subset of $A_{0}$.
\item The morphism $\varphi$ defines an isomorphism of monoïds in $(Set,\times,*)$, $A_{0}\backsimeq End_{A-mod}(A)$.
\item If there exists $f\in q_{0}$ which is invertible in $A_{0}$, then $q\backsimeq A$.  
\end{list}
\end{lemme}
Proof\\The proof of the three first items is clear. Let us prove the fourth.
The morphism $\varphi(f)$ is an automorphism of $A$ in $A-mod$ which factors through $j_{q}$. The morphism $j_{q}$ is thus both a monomorphism and a split epimorphism, hence an isomorphism.

\begin{definition}\label{def}
Let $\xymatrix{A\ar[r]^{u}&B}$ be a morphism of commutative monoids and $q$ be an ideal of $A$. 
The image of the ideal $q$ in $B$, denoted $qB$, is the image of the morphism $q\otimes_{A}B\rightarrow B$ in $B-mod$.
\end{definition}

\begin{lemme}\label{r1}
Let $\xymatrix{A\ar[r]^{u}&B}$ be a flat morphism of $Comm(\mathcal{C})$ (ie $-\otimes_{A}B$ is exact, cf \cite{8}), and $q$ be an ideal of $A$. Then $q\otimes_{A}B\backsimeq qB$.
\end{lemme}

\subsubsection{Multiplications}
\begin{definition}
Let $A$ be in $Comm(\mathcal{C})$, $f$ be an element of $A$ and $M$ an $A$-module.
The multiplication by $f$ in $M$, denoted $\tau_{f}$ is defined by the following commutative diagram
\begin{center}
$\xymatrix{M\otimes_{A}A\otimes k\ar[r]^{Id_{M}\otimes f}\ar[d]_{r_{M}\otimes Id_{k}}&M\otimes_{A}A\ar[d]^{r_{M}}\\M\otimes k\ar[r]_{\tau_{f}}&M}$
\end{center}
In particular, when $M=A$, $\tau_{f}=\varphi(f)$ where $\varphi$ is the adjunction morphism defined in the preliminaries.
\end{definition}

\begin{lemme}
Let $A$ be in $Comm(\mathcal{C})$, $M,N$ be two $A$-modules and $u$ be in $Hom_{\mathcal{C}}(M,N)$. The morphism $u$ is in $A-mod$ if and inly if it commutes with  morphisms $\tau_{f}$ for any element $f$ of a generating subset of $A$, i.e. verifies $\tau_{f}\circ (u\otimes Id_{k})=u\circ \tau_{f}$.
\end{lemme}
Proof\\It is easy to check, considering the facts that generating families of element are  epimorphic families, and functors $Z\otimes-$, $Z\in \mathcal{C}$ preserve epimorphisms.
\vskip 2pt
$\hfill\blacklozenge$

\begin{lemme}
Let $A$ be in $Comm(\mathcal{C})$, $q$ be an ideal of $A$. An element $f$ of $A$ is an element of $q$ if and only if $\varphi(f)$ factors through $j_{q}$.
\end{lemme}
Proof\\It comes from $\varphi(j_{q}\circ f)=j_{q}\circ \varphi(f)$.
\vskip 2pt
$\hfill\blacklozenge$ 

\subsection{Localizations}
\begin{definition}\label{d1}
Let $A$ be in $Comm(\mathcal{C})$ and $S$ be a multiplicative part of $A_{0}$. A localization of $A$ along $S$ is a couple $(S^{-1}A,w_{S})$, where $w_{S}:A\rightarrow S^{-1}A \in Comm(\mathcal{C})$, satisfying the following properties:
\begin{enumerate}
\item The image of the elements of $S$ are invertible, ie $(w_{S})_{0}(S)\subset (S^{-1}A)_{0}^{*}$. 
\item Any morphism verifying $i$ factors through $w_{S}$. More precisely, $\forall \xymatrix{A\ar[r]^{u}&B}\in Comm(\mathcal{C}),\;if\;u_{0}(S)\subset B_{0}^{*}$ then $\exists!\;v$ such that the following diagram commutes:
\begin{center}
$\xymatrix{A\ar[rr]^{u}\ar[rd]_{w_{S}}&&B\\&S^{-1}A\ar[ru]_{v}&}$
\end{center}
\end{enumerate}
If $S=<\{f\}>$ is the multiplicative part generated by one element $f$, $A_{f}$ will refer to $S^{-1}A$ and $w$ to $w_{S}$.
\end{definition}
\begin{remarque}\label{r4}
\begin{itemize}
\item Lemma \ref{r2} and proposition \ref{r3} prove the existence of $S^{-1}A$.
\item The morphism $A\rightarrow S^{-1}A$ is an epimorphism by construction.
\item It is well known that a morphism of monoids is an epimorphism if and only if the corresponding forgetfull functors on the categories of modules is fully faithfull. In particular, the category $S^{-1}A-mod$ is a full subcategory of $A-mod$.
\end{itemize}
\end{remarque}

\begin{proposition}\label{r3}
Let $A$ be a commutative monoid and $f$ be in $A_{0}$. Let $j_{1}$ and $j_{0}$ denote respectively the inclusion morphisms of $A=A^{\otimes_{A} 1}/S^{1}$ and $A=A^{\otimes_{A} 0}/S^{0}$ in $L_{A}(A)=\coprod_{n\in \mathcal{N}}A^{\otimes_{A} n}/S^{n}$ and define $\delta$ as $j_{1}\circ m_{A}\circ Id\otimes f\circ r_{A}^{-1}$. Then:
\begin{center}
$A_{f}\backsimeq coker(\xymatrix{L_{A}(A)\ar@<1ex>[r]^{\psi_{A}(\delta)}\ar@<-1ex>[r]_{\psi_{A}(j_{0})}&L_{A}(A)})$.
\end{center}
The colimit is taken in $Comm(\mathcal{C})$.
\end{proposition}
Proof :\\Let us prove first that this colimit satisfies $ii$ in definition \ref{d1}.\\Let $\xymatrix{A\ar[r]^{u}&B}$ be in $Comm(\mathcal{C})$ such that $g:=u_{0}(f)= u\circ f$ is invertible in $B_{0}$. Let $g^{-1}$ denote its inverse. Thus $m_{b}\circ (g\otimes g^{-1})\otimes r_{1}^{-1}= i_{B}$. Consider the morphism $L_{A}(A)\rightarrow B$ in $A-alg$ adjoint to $\beta = m_{b}\circ (u\otimes g^{-1})\circ r_{A}^{-1}:A\rightarrow B$. Then $\psi_{A}(\beta)\circ j_{0}$ is the unit morphism of $B$ in $A-alg$, ie $u$ and $\psi_{A}(\beta)\circ j_{1}=\beta$. Thus $\psi_{A}(\beta)\circ\delta = \beta\circ m_{A}\circ (Id\otimes f)\circ r_{A}^{-1}$. Moreover the following diagram is commutative:
\begin{center}
$\xymatrix{A\ar[d]_{r_{A}^{-1}}\ar@/^1pc/[rrd]^{\psi_{A}(\beta)\circ \delta}&&\\A\otimes 1\ar[r]^{Id_{A}\otimes r_{1}^{-1}}\ar[d]_{Id_{A}\otimes f}&A\otimes 1\otimes 1\ar[d]^{Id_{A}\otimes f\otimes Id_{1}}&b\\A\otimes A\ar[r]^{Id_{A}\otimes r_{A}^{-1}}\ar[d]_{m_{A}}&A\otimes A\otimes 1\ar[d]^{m_{A}\otimes Id_{1}}&B\otimes B\ar[u]_{m_{B}}\\A\ar[r]^{r_{A}^{-1}}&A\otimes 1\ar[r]_{Id_{A}}&A\otimes 1\ar[u]_{u\otimes g^{-1}}}$
\end{center}
Then using compatibility between $\otimes$ and $\circ$, commutative monoid morphism's structure of $u$ and associativity of multiplication, the computation of $\beta\circ m_{A}\circ (Id\otimes f)\circ r_{A}^{-1}$ gives $u$. Thus $\psi_{A}(\beta)\circ\delta=\psi_{A}(\beta)\circ j_{0}$, hence $\psi_{A}(\beta)\circ\psi_{A}(\delta)=\psi_{A}(\beta)\circ \psi_{A}(j_{0})$. Finally, by the universal property of colimit, there exists a unique morphism from this colimit to $B$.\vskip 3pt
Let us prove now that the colimit satisfies $i$ in definition \ref{d1}. We denote now this colimit by $A_{f}$, and the morphism $A\rightarrow A_{f}$ by $\omega$. There exists a morphism $\bar{p}$ from $L_{A}(A)$ to $A_{f}$ such that $\bar{p}\circ \delta = \omega$. Let $p:A\rightarrow A_{f}$ denote the adjoint of $\bar{p}$. Then $p\circ m_{A}\circ (Id\otimes f)\circ r_{A}^{-1}=p$. Moreover, considering the identity $m_{A_{f}}\circ (\omega\otimes (p\circ i_{A}))\circ r_{A}^{-1}= p$, it is clear that the morphism $p$ plays the same role as $\beta$, and  $p\circ i_{A}$ will play the role of $g^{-1}$. Thus, an analogous computation gives : 
\begin{center}
$\omega\circ i_{A} = p\circ m_{A}\circ (Id\otimes f)\circ r_{A}^{-1}\circ i_{A} = m_{A_{f}}\circ ((\omega\circ f)\otimes p\circ i_{A})\otimes r_{1}^{-1}$.
\end{center}
Thus the product of $\omega\circ f$ and $p\circ i_{A}$ is the unit element of $(A_{f})_{0}$.

\begin{lemme}\label{r2}
A category, still denoted $S$, is associated to $S$. Its elements are elements of $S$ and for two elements $f,g$ of $S$, the morphisms from $f$ to $g$ are the elements $h$ of $S$ verifying $f*h=g$. There is an isomorphism $S^{-1}A\backsimeq Colim_{S}A_{f}$. Furthermore, as $S$ is multiplicative, this colimit is filtered. 
\end{lemme}
Proof\\By the universal property of objects $S^{-1}A$ and $A_{f}$.
\vskip 2pt
$\hfill\blacklozenge$

\begin{proposition}\label{colfil}
Let $A$ be a commutative monoid and $f$ be in $A_{0}$.\\The object $B:= colim(\xymatrix{A\ar[r]^{\varphi(f)}&A\ar[r]^{\varphi(f)}&A\ar[r]^{\varphi(f)}&...})$ can be provided with a monoid structure and is isomorphic in $A-alg$ to $A_{f}$.
\end{proposition}  
Proof:\\Let $B$ denote this colimit. First, $B$ can be provided with a monoid structure. The unit morphism $1\rightarrow B$ is induced by the morphism from $1$ to the first term $A$ in the colimit. Let $\mathcal{N}^{+1}$ denote the category whose objects are positive integers and whose morphisms are couples $(i,i+1)$. Let $F$ denote the functor from $\mathcal{N}^{+1}$ to $A-mod$ which sends any integer to $A$ and any morphim to $\varphi_{f}$. Then, clearly,  $B\otimes_{A}B\backsimeq Colim (F\otimes_{A}F)\backsimeq Colim F\backsimeq B$. This isomorphism defines the multiplication morphism of $B$. One can check easily that these morphisms provide $B$ with a structure of monoid in $\mathcal{C}$.\vskip 6pt

Next step is to prove that $g:=i_{B}\circ f$ is invertible in $Hom_{A-mod}(A,B)\backsimeq B_{0}$. There exists a morphism from the second term $A$ in the colimit to $B$, denoted arbitrarily $g^{-1}$. It verifies $g^{-1}\circ \varphi(f)=i_{B}$. As $\varphi(f)=m_{A}\circ(Id_{A}\otimes_{A}f)\circ r_{A}^{-1}$, we get
\begin{center}
$i_{B}=g^{-1}\circ \varphi(f)=g^{-1}\circ m_{A}\circ (Id_{A}\otimes_{A}f)\circ r_{A}^{-1}$
\end{center}
and by $A$-linearity of $g^{-1}$:
\begin{center}
$=\mu_{B}\circ(g^{-1}\otimes_{A}Id_{A})\circ (Id_{A}\otimes_{A}f)\circ r_{A}^{-1}$
\end{center}
where $\mu_{B}=m_{B}\circ (Id_{B}\otimes_{A}i_{B})$.\vskip 6pt
Thus
\begin{center}
$=\mu_{B}\circ (g^{-1}\otimes_{A} f)\circ r_{A}^{-1}$
\end{center}
so
\begin{center}
$m_{B}\circ (g^{-1}\otimes_{A} g)\circ r_{A}^{-1}$.
\end{center}
\vskip 6pt
Last but not least, let us prove that $B$ satisfies property $ii$ of the definition of $A_{f}$.
Let $v$ denote the unit morphism of $B$ in $A-alg$. Let $h:A\rightarrow C$ be a morphism in $Comm(\mathcal{C})$ such that $g:=h_{0}(f)$ is invertible in $C_{0}$. Let $g^{-1}$ denotes the inverse of $g:=h\circ f$ in $C_{0}$. Then the following diagram commutes:
\begin{center}
$\xymatrix{A\ar[r]^{r_{A}^{-1}}&A\otimes 1\ar[r]^{h\otimes g^{-1}}&C\otimes C\ar[r]^{m_{C}}&C\\A\ar[u]^{\varphi(f)}\ar[rrru]_{h}}$
\end{center}
Indeed $m_{C}\circ (h\otimes g^{-1})\circ r_{A}^{-1}$ is just $\varphi(g^{-1})$ and $h\circ \varphi(f)$ is $\varphi(g)$. So this diagram gives their product in $Hom_{A-mod}(A,C)$ which is the unit morphism of $C$ in $A-alg$, $h$. Using the same method with inverses $g^{\otimes-n}$ of $g^{\otimes n}$, we contruct a family of morphisms from $A$ to $C$ which makes the diagram of the colimit $B$ commutes. Thus there is an induced morphism in $A-mod$, $s:B\rightarrow C$ such that $s\circ v=h$.\\The morphism $s$ is in fact a morphism of $A-alg$. Indeed, the commutativity of the diagram:
\begin{center}
$\xymatrix{B\ar[r]^{s}&C\\B\otimes_{A}B\ar[u]^{m_{B}}\ar[r]_{s\otimes_{A}s}&C\otimes_{A}C\ar[u]_{m_{C}}}$
\end{center}
is equivalent, by the universal property of colimit, to the commutativity, for each $n$, of the diagram:
\begin{center}
$\xymatrix{&B\ar[r]^{s}&C\\&B\otimes_{A}A\ar[u]_{\mu_{B}\circ(\varphi(f)^{n}\otimes_{A}Id_{A})}\ar[r]_{s\otimes_{A}\phi(g^{-n})}&C\otimes_{A}C\ar[u]_{m_{C}}\\B\ar[ur]_{r_{B}^{-1}}\ar@/^-1pc/[urr]\ar@/^1pc/[uur]^{\tau_{g^{-n}}}&&}$
\end{center}
which is clear.

\begin{corollaire}\label{r5}
Let $q$ be an $A-module$, $f$ be in $A_{0}$, and $S$ be a multiplicative part of $A_{0}$, then 
\begin{enumerate}
\item $q_{f}:=q\otimes_{A}A_{f}\backsimeq colim(\xymatrix{q\ar[r]^{\tau(f)}& ...})$.
\item $q\otimes_{A}S^{-1}A\backsimeq colim_{S}q_{f}$.
\end{enumerate}
\end{corollaire}
Proof\\It is implied by lemma \ref{r2}, proposition \ref{colfil} and the fact that $q\otimes_{A}-$ commutes with colimits.
\begin{corollaire}
Let $S$ be a multiplicative part of $A_{0}$. The category $S^{-1}A-mod$ of modules over $S^{-1}A$ is exactly the subcategory of $A-mod$ of modules $M$ such that for any element $f$ of $S$, $\tau_{f}$ is an invertible endomorphism of $M$. 
\end{corollaire}
Proof:\\ Clearly, if $q$ is in $S^{-1}A-mod$, it verifies the property. Let $q$ be an $A$-module verifiying the property. Then, by previous corollary $q_{f}\backsimeq q$ $\forall f\in S$ hence $q\otimes_{A}S^{-1}A\backsimeq q$ and thus $q$ is in $S^{-1}A-mod$.
\vskip 6pt
The localizations of $A$ are important because they are the basic examples of Zariski open subobjects of $Spec(A)$. Thus, next step is to prove it.

\subsection{Zariski open objects}
\begin{definition}(cf \cite{8})
\begin{enumerate}
\item A morphism $A\rightarrow B$ in $Comm(\mathcal{C})$ is $\mathcal{C}$-flat if the functor $-\otimes_{A}B$ commutes with finite enriched limits (cf definition \ref{a2} in appendix).
\item A morphism $A\rightarrow B$ in $Comm(\mathcal{C})$ is finitely presented if $B$ is finitely presented in $A-alg$.
\item A morphism $A\rightarrow B$ in $Comm(\mathcal{C})$ is a formal Zariski open morphism if it is both a $\mathcal{C}$-flat morphism and an epimorphism.
\item A morphism $A\rightarrow B$ in $Comm(\mathcal{C})$ is a Zariski open morphism iff it is a finitely presented formal Zariski open morphism. In this case, $Spec(B)$ will be called a Zariski open subobject of $Spec(A)$.
\end{enumerate}
\end{definition}
\begin{remarque}
This definition is not exactly the definition of \cite{8}. Indeed, $\mathcal{C}$-flat means that the functor is not only exact but also commutes with finitely presented exponentiation (see proposition \ref{a3} in the appendix). This more precise definition is necessary to prove the theorem and is also a generalisation of usual cases. This hypothesis is also necessary to recover the notion of differential scheme (no comparison had been writen for the moment.). In facts, there is an equivalence between flat and $\mathcal{C}$-flat if finitely presented objects in $A-mod$ are quotients of free objects of finite type. In this case, exponential objects depend only on finite limits. One example is the case of abelian groups $(Z-mod,\otimes,Z)$ in which for any ring $A$ and finitely presented $A$-module $M$, there exists a resolution $A^{q}\rightarrow A^{p}\rightarrow M$. 
\end{remarque}
\begin{lemme}\label{r6}Let $A$ be in $Comm(\mathcal{C})$.
\begin{enumerate}
\item If $X\in \mathcal{C}$ (resp $A-mod$) is finitely presented then so is $L_{A}(X)\in Comm(\mathcal{C})$ (resp $A-alg$).
\item Finitely presented objects are stable under finite colimits. 
\end{enumerate}
\end{lemme}
Moreover, in $A-mod$, filtered colimits commute with finite enriched limits (Corollary \ref{a1} in appendix). In particular, localizations are Zariski open objects. Indeed, the following corollary results:    

\begin{corollaire}(of \ref{r3}, \ref{r4}, \ref{r5} and \ref{r6})\\
Let $A$ be in $Comm(\mathcal{C})$, $f$ be in $A_{0}$, and $S$ be a multiplicative part of $A_{0}$, then $A\rightarrow S^{-1}A$ is a formal Zariski open morphism and $A\rightarrow A_{f}$ is a Zariski open morphism, i.e. $spec(S^{-1}A)$ is a formal Zariski open subobject of $Spec(A)$ and $spec(A_{f})$ is a Zariski open subobject of $Spec(A)$. 
\end{corollaire}

\section{Gabriel Filters}
 
\subsection{Gabriel filters and Borceux-Quinteiro's theorem}
The key idea of the classification is to indroduce the notion of Gabriel filter to make a link between Zariski open objects and ideals. A Gabriel filter on a monoid $A$ will be a topology of Grothendieck on the category $\mathcal{B}A$, defined here in an enriched sense.  
\begin{definition}Let $A$ be in $Comm(\mathcal{C})$
\begin{enumerate}
\item Let $\mathcal{B}A$ denote the $\mathcal{C}$-category consisting of one object $*$ and an object of morphisms $\underline{Hom}_{\mathcal{B}A}(*,*):=A\in \mathcal{C}$.
\item Let $Pr(\mathcal{B}A)$ denote the $\mathcal{C}$-category of $\mathcal{C}$-presheaves i.e. of $\mathcal{C}$-functors from $\mathcal{B}A$ to $\mathcal{C}$.
\item The Yoneda $\mathcal{C}$-functor consists of the $\mathcal{C}$-functor $*\rightarrow h_{*}$ which is defined on objects by $h_{*}(*)=A$ and on morphisms by the adjoint $h_{*} :A\rightarrow \underline{End}_{\mathcal{C}}(A)$ of $m_{A}$.
\end{enumerate}
\end{definition}
In fact, any $\mathcal{C}$-presheaf on $\mathcal{B}A$ sends $*$ to an $A-module$ and the functor
\begin{center}
$F:Pr(\mathcal{B}A)\rightarrow A-mod$ \vskip 6pt $(G:\mathcal{B}A\rightarrow \mathcal{C})\rightarrow G(*)$
\end{center}
defines a $\mathcal{C}$-equivalence from $Pr(\mathcal{B}A)$ to $A-mod$.
\begin{definition}(cf \cite{2} for general definition)
A Grothendieck's $\mathcal{C}$-topology on $\mathcal{B}A$ consists of a set $J$ of ideals of $A$ such that:
\begin{enumerate}
\item $A\in J$
\item $\forall q\in J$,  $\forall k\in \mathcal{C}_{0}$ and  $\forall f\in A_{k}$, $f^{-1}(q)=(q^{k}\times_{A^{k}}A)\in J$.
\item Let $q,r$ be two ideals of $A$, $q\in J$. If $\forall f:k\rightarrow q,\;k\in \mathcal{C}_{0}$ the ideal $f^{-1}(r)$ is in $J$, then so is $r$. 
\end{enumerate}
Where $q^{k}=\underline{Hom}(k,q)$ and $A^{k}=\underline{Hom}(k,A)$ are expoential objects.\vskip 2pt This topology will be called a Gabriel filter. In $ii$ the morphism from $A$ to $\underline{Hom}_{\mathcal{C}}(k,A)$ is the adjoint of the morphism $\varphi(f):(A\otimes k)\rightarrow A$.
\end{definition}

\begin{lemme}\label{filii}
Let $G$ be a Gabriel filter, $q$ be an ideal in $G$, $k$ be in $\mathcal{C}_{0}$ and  $f$ be in $A_{k}$. The ideal $f^{-1}q$ is the $sup$ (or union), in the poset of ideals of $A$ containing $q$,  of the ideals $p$ such that the multiplication $\tau_{f}:p\rightarrow p$ factors through $q\hookrightarrow p$. This family of ideals is denoted $cont_{f}q$.
\end{lemme}
Proof\\
Weak underlying set functors commute with limits and the following commutative diagram proves that $f^{-1}q\in cont_{f}q$
\begin{center}
$\xymatrix{A_{k'}\ar[r]^{\tau_{f}^{*}}&A_{k'\otimes k}\\f^{-1}q\ar@{^{(}->}[u]\ar[r]&q_{k'\otimes k\ar@{^{(}->}[u]}}$
\end{center}
Moreover, if $q'$ is in $cont_{f}q$, it is easy to check , by the universal property of $f^{-1}q$, that $q'\hookrightarrow f^{-1}q$.
\vskip 2pt
$\hfill\blacklozenge$ 

\begin{definition}
A localization of $A-mod$ consists of a full reflexive (saturated) sub-$\mathcal{C}$-category of $A-mod$ such that the left adjoint of the inclusion $\mathcal{C}$-functor is left $\mathcal{C}$-exact.
\end{definition}
\begin{theoreme}
(Borceux-Quinteiro)\\Let $A$ be in $Comm(\mathcal{C})$, there is a bijection between the set of Gabriel filters of $A$ and the set of localizations of $A-mod$.
\end{theoreme}
The poset of Gabriel filters of $A$ will be denoted $GabA$, the poset of localizations of $A-mod$ will be denoted $LocA$. 
\begin{lemme}\label{corfl} (cf \cite{2})
\begin{list}{$\triangleright$}{}
\item Let $\mathcal{L}$ be a localisation of $A-mod$. The Gabriel filter associated to $\mathcal{L}$, denoted $G_{\mathcal{L}}$ consists of the set of ideals $q$ of $A$ verifying for all $M$ in $\mathcal{L}$
\begin{center}
$\underline{Hom}_{A-mod}(A,M)\backsimeq \underline{Hom}_{A-mod}(q,M)$
\end{center}
\item Let $G$ be a Gabriel filter. The localisation associated to $G$, denoted $\mathcal{L}_{G}$, is the full sub-$\mathcal{C}$-category of $A-mod$ whose objetcs are $A$-modules $M$ verifying for all $q$ in $G$
\begin{center}
$\underline{Hom}_{A-mod}(A,M)\backsimeq \underline{Hom}_{A-mod}(q,M)$
\end{center}
\end{list}
\end{lemme}

\begin{remarque}\label{r7}
\begin{list}{-}{}
\item The bijection $GabA\backsimeq LocA$ is a contravarient isomorphism of posets. 
\item The poset $LocA$ has small intersections. If $\mathcal{L}_{i}$ is a small family of localisations with left adjoint functors $l_{i}$, its intersection consists of the sub-$\mathcal{C}$-category of $A-mod$ whose set of objects is $\cap_{i}Ob(\mathcal{L}_{i})$ and whose left adjoint functor is $Colim(l_{i})$. 
\item The poset $GabA$ has small intersections. A Small intersection in $GabA$ is an intersection in the sense of sets. 
\item The contravarient isomorphism of poset $GabA\backsimeq LocA$ implies that $LocA$ and $GabA$ have small unions.  
\end{list}
\end{remarque}

\begin{proposition}
Let $Y:=Spec(B)$ be a formal Zariski open subobject of $X:=Spec(A)\in Aff_{\mathcal{C}}$. The category $B-mod$ is a localization of $A-mod$. Its associated Gabriel filter is the set of ideals $q$ such that $q\otimes_{A}B\backsimeq B$.
\end{proposition}
Proof\\
The category $B-mod$ is clearly a localization of $A-mod$. The caracterization of \ref{corfl}, induces by adjunction, for all $M\in B-mod$, an isomorphism
\begin{center}
$\underline{Hom}_{B-mod}(B,M)\backsimeq \underline{Hom}_{B-mod}(q\otimes_{A}B,M)$.
\end{center}
Thus, by Yoneda $q\otimes_{A}B\backsimeq B$.
\vskip 2pt
$\hfill\blacklozenge$ 
\vskip 2pt
The object $Spec(S^{-1}A)$ is a formal Zariski open subobject of $Spec(A)$. Thus, there is an associated localization of $\underline{A-mod}$. The following proposition describes its associated Gabriel filter.
\begin{proposition}\label{r9}
Let $A$ be a Commutative monoid, $f\in A_{0}$ and S a multiplicative part of $A_{0}$.
\begin{enumerate}
\item The Gabriel filter $G_{f}$ associated to $A_{f}$ is the set of ideals of $A$ containing a power of $f$.
\item The Gabriel filter $G_{S}$ associated to $S^{-1}A$ is $\cup_{f\in S}G_{f}$
\end{enumerate}
\end{proposition}
Proof\\Let $G$ denote the set of ideals of $A$ containing a power of $f$. Let $q$ be an element of $G$. There exists $n$ such that $f^{n}\in q_{0}$, thus $q_{f}$ contains an invertible element of $A_{f}$ and by $\ref{r0}$ is isomorphic to $A_{f}$. Therefore, $G\subset G_{f}$.\vskip 4pt
Reciprocally, let $q$ be in $G_{f}$, then $i_{q_{f}}:q_{f}\backsimeq q\otimes_{A}A_{f}\rightarrow A_{f}$ is an isomorphism and the following diagram is commutative:
\begin{center}
$\xymatrix{A_{0}\ar[rr]&&(A_{f})_{0}\backsimeq (A_{0})_{f}\\q_{0}\ar@{^{(}->}[u]^{-\circ i_{q}}\ar[rr]&&(q_{f})_{0}\backsimeq (q_{0})_{-\circ\tau_{f}}\ar[u]}$
\end{center}
As this diagram is in Sets, $\exists n\in \mathcal{N}$ such that $f^{\otimes n}\in q_{0}$. Thus $q\in G$. Finally $G=G_{f}$.\\For $S^{-1}A$, clearly $\cup_{f\in S}G_{f}\subset G_{S}$ and the other inclusion is obtained with an analogous commutative diagram in Set.

\subsection{Primitive Quasi-Compact and Locally Primitive Gabriel Filters}

\begin{definition}Let $A$ be a commutative monoid 
\begin{list}{$\triangleright$}{}
\item Let $q,q'$ be two ideals of $A$. The product $q.q'$ is the image of the morphism $q\otimes q'\rightarrow A$.
\item An ideal $p$ of $A$ is prime if for all $q,q'$ sucht that $q.q'\subset p$, $q\subset p$ or $q'\subset p$.
\item A set $G$ of ideal of $A$ is filtered if for all $q\subset q'$, $q\in G\Rightarrow q'\in G$.  
\item A set of ideals of $A$ is commutative if it is closed under product of ideals.
\item A set of ideals $G$ is \textit{quasi-compact} if, for any ideal $q$ in $G$ isomorphic to a filtered colimit $colim_{i\in \mathcal{I}}q_{i}$, there exists $i$ such that $q_{i}$ is in $G$.  
\end{list}
\end{definition}
\begin{lemme}\label{r11}Let $A$ be a commutative monoid
\begin{list}{$\diamond$}{}
\item A commutative Gabriel filter is a monoid in Sets.
\item Any Gabriel filter is filtered.
\item Any Gabriel filter is commutative.
\item A finite union or a finite intersection of quasi-compact sets of ideals is quasi-compact.
\item A Gabriel filter $G$ is quasi-compact if and only if for any ideal $q$ in $G$ there exists an ideal $q'\subset q$, finitely generated, in $G$.
\item Let $A$ be a commutative monoid, $q,q'$ be two ideals of $A$ and $F$ be a generating family of $q$. There is an isomorphism $q.q'\backsimeq \cup_{f\in F}Im(\tau_{f})$ where $\tau_{f}:q'\rightarrow q'$. 
\end{list}
\end{lemme}
Proof\\The first property is clear. The second property comes from axiom $ii$ in the definition of Gabriel filter. The third comes from axiom $iii$. The fourth property is due to the fact that finite unions and intersections of Gabriel filters are unions and intersections in the sense of sets.\vskip 2pt
Let us now prove the fifth property. If $G$ is a quasi-compact Gabriel filter, any ideal $q$ of $G$ is the filtered colimit of its finitely generated sub-ideals (\ref{prgen}) thus there is a finitely generated sub-ideal of $q$ in $G$. Reciprocally, if $q\backsimeq colim_{\mathcal{I}}(q_{i})$ is a filtered colimit in $G$ along $\mathcal{I}$ and $q'$ is a finitely generated sub-ideal of $q$ in $G$. The ideals $q'$ is of finite type and the morphism from $q'$ to the colimit factors through an ideal $q_{i}$. The morphism $q\rightarrow q_{i}$ is clearly a monomorphism, thus $q_{i}\in G$.\vskip 2pt
Let us prove now the sixth property. There is an isomorphism $q\backsimeq Im(K(F)\rightarrow q)$. Moreover the image of $K(F)$ by this morphism is a colimit hence commutes with tensor product. Furthermore $q'\backsimeq K(i(q')\backsimeq Im(K(i(q'))\rightarrow q'))$ where $i$ is the Yoneda functor from $\mathcal{C}$ to $Pr(\mathcal{C}_{0})$. By \ref{prgen}, the image of $K(F)\otimes q'\rightarrow q\otimes q'$ is $q\otimes q'$.\vskip 4pt Thus
\begin{center}
$q\otimes q'\backsimeq Im(Colim_{(k,f),(k',f')\in \mathcal{C}^{F}_{0}\times \mathcal{C}_{0}^{i(q')}}k\otimes k'\rightarrow q\otimes q')$.
\end{center}
In particular, $F\times i(q')$ generates $q.q'$. Moreover for all $k\in \mathcal{C}_{0}$
\begin{center}
$k\otimes q'\backsimeq Colim_{(k',f')\in \mathcal{C}^{i(q')}_{0}}k\otimes k'$.
\end{center}
Thus $Im(\tau_{f})$ is generated by $\{f\}\times i(q')$. Therefore, any generator of $q.q'$ is included in an object $Im(\tau_{f})$, $f\in F$ and finally $\cup_{f\in F}Im(\tau_{f})\backsimeq q.q'$.
\vskip 2pt
$\hfill\blacklozenge$

The elements of an ideal play a particularly important role. We describe them more precisely.

\begin{definition}Let $A$ be a commutative monoid
\begin{list}{$\triangleright$}{}
\item An isomorphism between two elements $f\in A_{k},\;g\in A_{k'}$ is an isomorphism from $k$ to $k'$ equalizing the two elements.
\item Let $f$ be an element of $A$, the ideal generated by $f$ is denoted $(f)$.
\end{list}
\end{definition}
\begin{lemme}
Let $A$ be a commutative monoid. The collection of isomorphism classes of elements of $A$ is a symmetric monoid in sets, its unit is $i_{A}$. The multiplication is denoted $f\oslash g$. Moreover $(f\oslash g)=(f).(g)$. In particular $(f^{\oslash n})=(f)^{n}$. 
\end{lemme}
Proof\\The product given by $f,g\rightarrow f\oslash g:=f\circ (g\otimes Id_{k})$ provides the set of isomorphism classes of elements with a structure of commutative monoid. The formula $(f\oslash g)=(f).(g)$ is a particular case of the last property proved in previous lemma. 
\vskip 6pt
We will now give a more convenient equivalent definition for these commutative Gabriel filters, in particular those of finite type. 

\begin{proposition}\label{o1}Let $A$ ba a commutative monoid
\begin{list}{$\triangleright$}{}
\item A Gabriel filter of $A$ consists of a subset $G$ of the set of ideals of $A$ such that
\begin{enumerate}
\item The set $G$ is non empty.
\item The set $G$ is filtered.
\item The property $iii$ of the definition of Gabriel filter holds in $G$.
\end{enumerate}
\item A Gabriel filter of finite type of $A$ consists of a subset of the set of ideals of $A$ such that:
\begin{enumerate}
\item The set $G$ is non empty.
\item The set $G$ is filtered.
\item The set $G$ commutative.
\item The set $G$ is of finite type. 
\end{enumerate}
\end{list}
\end{proposition}
Proof:\\Any Gabriel filter is non empty and filtered. Reciprocally, we need to check that property $i$ and $ii$ hold in $G$. $G$ is filtered and non empty thus contains $A$. Let $q$ be in $G$ and $f\in A_{k}$, the ideal $f^{-1}q$ contains $q$ thus is in $G$.\vskip 6pt
Let us now prove the caracterisation of Gabriel filters of finite type. Let $G$ be a set of ideals verifying these four properties. We have already proved that properties $i$ and $ii$ of the definition of Gabriel filter hold in $G$, let us prove now it for $iii$. Let $q$ be in $G$ and $q'$ be such that $\forall f\in_{k} q$, $f^{-1}q'$ is in $G$. Then let $I$ be the family of finite subets of elements of $q$, i.e. of $\prod_{k\in \mathcal{C}_{0}}Hom_{\mathcal{C}}(k,q)$, let $q_{X}$ denotes the ideal generated by a subset $X$ of $\prod_{k\in \mathcal{C}_{0}}Hom_{\mathcal{C}}(k,q)$. By \ref{prgen}, $q\backsimeq colim_{i\in I}q_{i}$. As $G$ is finitely presented, there exists $X$ such that $q_{X}$ is in $G$.\\The following diagram is commutative:
\begin{center}
$\xymatrix{A\otimes k \ar[r]^(0.6){\tau_{f}}& A\\f^{-1}(q')\otimes k \ar[u]\ar[r]_(0.6){\tau_{f}}& q'\ar[u]_{i_{q'}}}$
\end{center}
Thus, $q_{X}.\cap_{f\in X}f^{-1}q'\backsimeq \cup_{f\in X}Im_{\cap_{f\in X}f^{-1}q'}(\tau_{f})\subset q'$. As $G$ is closed under finite intersections and finite products, $q'$ is in $G$.
\vskip 2pt
$\hfill\blacklozenge$

\begin{proposition}
Let $u:A\rightarrow B$ be in $Comm(\mathcal{C})$. If $u$ is a formal Zariski open morphism then the Gabriel filter associated to $B$ is quasi-compact.
\end{proposition}
Proof\\Let $q$ be in $G$ and $\mathcal{I}$ be a filtered diagram such that $q\backsimeq colim_{i\in \mathcal{I}}q_{i}$. There is an isomorphism $q\otimes_{A}B\backsimeq B$ hence $B\backsimeq colim_{i\in \mathcal{I}}(q_{i}\otimes_{A}B)$. As $B$ is finitely presented in $A-mod$, there exists $i$ such that this morphism factors through $q_{i}\otimes_{A}B$. As $u$ is flat, $q_{i}\otimes_{A}B\backsimeq q_{i}.B$ is an ideal of $B$. Finally, by lemma \ref{r0} $q_{i}\otimes_{A}B\backsimeq B$, hence $q_{i}\in G$.
\vskip 2pt
$\hfill\blacklozenge$

\begin{definition}
Let $A$ be a commutative monoid and $q$ be an ideal of $A$. 
\begin{list}{$\triangleright$}{}
\item The Gabriel filter $G_{q}$ is the smallest Gabriel filter of $A$ containing $q$ in the poset $GabA$. A Gabriel filter generated by one ideal in this sense will be called a primitive Gabriel filter. 
\item The Gabriel filter $G^{q}$ is defined by $G^{q}=\cap_{f\in q}G_{f}$ in the poset $GabA$. A Gabriel filter generated by one ideal in this sense will be called a locally primitive Gabriel filter. 
\end{list}
\end{definition}
\begin{remarque}
The set of prime ideals included in $G^{q}$ is the set of prime ideals containing $q$.
\end{remarque}

\begin{lemme}\label{r12}
Let $A$ be a commutative monoid and $q$ be an ideal of $A$ of finite type. The Gabriel filter $G_{q}$ consists of the set of ideals of $A$ containing a power of $q$.
\end{lemme}
Proof\\Let $G$ be the set of ideals containing a power of $q$. By \ref{o1}, it is a quasi-compact Gabriel filter thus $G_{q}\subset G$. Moreover its elements are obtained from $q$ by products and inclusions, therefore $G\subset G_{q}$. Finally $G=G_{q}$.
\vskip 2pt
$\hfill\blacklozenge$

\begin{corollaire}\label{r13}
Let $A$ be a commutative monoid and $q$ an ideal of $A$
\begin{list}{$\diamond$}{} 
\item If $F$ is a generating family of $q$ then $G^{q}=\cap_{f\in F}G_{f}$.  
\item If $q=\cup_{i\in I}q_{i}$, then $G^{q}=\cap_{i\in I}G^{q_{i}}$. 
\end{list}
\end{corollaire}
Proof\\ The proofs of these two properties are equivalent so let us prove the first. There is a clear inclusion $G^{q}\subset \cap_{f\in F}G_{f}$. Let $p$ be an ideal in $\cap_{f\in F}G_{f}$.
Let $g$ be an element of $q$. Recall that $(g)$ refers to the ideal generated by $g$. By \ref{prgen}, there is an isomorphism
\begin{center}
$q\backsimeq colim_{F'\subset F,\;fini}(\cup_{f\in F'} q_{F'})$
\end{center}
As $(g)$ is an ideal of finite type , the monomorphism $(g)\rightarrow q$ factors through an object $q_{F'}$ in the colimit. There is an inclusion $\cap_{f\in F}G_{f}\subset G^{q_{F'}}$, thus $p$ is in $G^{q_{F'}}$. As $F'$ is finite $G^{q_{F'}}=G_{q_{F'}}$ then $p$ contains a power of $q_{F'}$. As $g$ is an element of $q_{F'}$, $p$ contains a power of $g$. This is true for all $g\in q$, thus $G^{q}=\cap_{f\in F}G_{f}$.
\vskip 2pt
$\hfill\blacklozenge$ 

\begin{lemme}\label{prime}
Let $G$ and $G'$ be two locally prmitive (resp primitive) Gabriel filters generated respectively by $q$ and $q'$. Their union in the poset of locally primitive (resp primitive) Gabriel filters is the locally primitive (resp primitive) Gabriel filter generated by $q.q'$. In particular, the prime ideals contained in the union are exactly the prime ideals contained in one of the two Gabriel filters.
\end{lemme} 
Proof\\In the primitive case, $G_{q.q'}$ contains $G$ and $G'$ and for any ideal $p$ such that $G_{p}$ contains $G$ and $G'$, $G_{q.q'}$ is clearly included in $G_{p}$.\vskip 2pt In the locally primitive case, we write $G^{q}\cup G^{q'}=\cap_{f\in q,f'\in q'}G_{f}\cup G_{f'}=G^{q}\cup G^{q'}=\cap_{f\in q,f'\in q'}G_{f\oslash f'}=G^{q.q'}$.
\vskip 2pt
$\hfill\blacklozenge$ 

\section{Zariski Open Objects Classififcation}

\subsection{Zariski Open subobjects}

\subsubsection{Affine Zariski Open subobjects}

\begin{lemme}\label{o3}
Let $A$ be a commutaive monoid and $spec(B),spec(B')$ be two Zariski open subobjects of $Spec(A)$. Then $spec(B\otimes_{A}B')$ is a Zariski open subobject of $Spec(A)$ and $B\otimes_{A}B'-mod=B-mod\cap B'-mod$. The intersection is taken in the poset $LocA$.
\end{lemme}
Proof\\ The localizations $B-mod\cap B'-mod$ and $B\otimes_{A}B'-mod$ have clearly the same objects hence are equivalent $\mathcal{C}$-categories. The two left adjoint functors are thus adjoint to the same forgetfull functor, thus are isomorphic.
\vskip 2pt
$\hfill\blacklozenge$

We give now a first interresting result which make a link between Gabriel filters and Zariski open subobjects.

\begin{theoreme}\label{CTF}Finite covering theorem\\
Let $A$ be a commutative monoid and $(Spec(B_{i}))_{i\in I}$ be a finite family of Zariski open subobjects of $Spec(A)$. This family is a covering of the Zariski open subobject $spec(B)$ for the topology defined on $Comm(\mathcal{C})$ (cf \cite{8}) if and only if $B-mod=\cup_{i\in I}(B_{i}-mod)$, where the union is taken in the poset $LocA$. 
\end{theoreme}
We define a functor $l:A-mod\rightarrow A-mod$ 
\begin{center}
$l(M):=ker(\xymatrix{\prod_{i,j}M\otimes_{A}B_{i}\otimes_{B}B_{j}\ar@<1ex>[r]\ar@<-1ex>[r]&\prod_{i}M\otimes_{A}B_{i}})$
\end{center}
This functor commutes with finite enriched limits (all the monoids in the kernel are $\mathcal{C}$-flat). It also verifies, for all $M\in A-mod$ and all $i$, $l(M)\otimes_{A}B_{i}\backsimeq l(M\otimes_{A}B_{i})$ by flatness of $B_{i}$.\vskip 2pt
Let us prove first that $l$ is isomorphic to identity on the categorie $B_{i}-mod$ for all $i$. Let $M$ be in $B_{i}-mod$, there is a unique monomorphism $M\rightarrow l(M)$ equilizing the two arrows of the kernel, induced by the natural morphisms $M\rightarrow M\otimes_{A}B_{i}\otimes_{A}B_{j}$. Furthermore
\begin{center}
$l(M)\rightarrow \prod_{j\in I}M\otimes_{A}B_{j}\rightarrow M\otimes_{A}B_{i}\backsimeq M$ 
\end{center}
Thus the composition on $M$ is the identity and the composition on $l(M)$ is the unique endomorphism of $l(M)$ equalizing the two arrows of the kernel hence is identity. Finally $M\backsimeq l(M)$.\vskip 2pt Let us now verify that $l²=l$. For all $M\in A-mod$
\begin{center}
$l²(M)\backsimeq Ker(\xymatrix{\prod_{i,j}l(M)\otimes_{A}B_{i}\otimes_{A}B_{j}\ar@<1ex>[r]\ar@<-1ex>[r]&\prod_{i}l(M)\otimes_{A}B_{i}})$
\end{center}
And for all $i$ 
\begin{center}
$l(M)\otimes_{A}B_{i}\backsimeq l(M\otimes_{A}B_{i})\backsimeq M\otimes_{A}B_{i}$
\end{center}
Thus $l²(M)$ is isomorphic to $l(M)$.\vskip 4pt
We prove then that the category $\mathcal{L}$ defined as the essential image of $l$ is a localization. First, we prove that $l$ is the left adjoint  of the forgetfull functor. The unit of the adjunction is
\begin{center}
$M\rightarrow l(M)$
\end{center}
The counit of the adjunction, for $N\backsimeq l(M)$ in the essential image of $l$, is the isomorphism 
\begin{center}
$l(N)\backsimeq N$
\end{center}
The triangular identities for the unit and the counit are obtained from the triangular identities of the unit and the counit of the adjunctions between forgetfull functors and functors $-\otimes_{A}B_{i}$. 
\vskip 3pt
Finally, $\mathcal{L}$ is a localization with left adjoint functor $l$. It contains all the $B_{i}-mod$ and thus their union in $LocA$. It is by definition contained in $B-mod$. 
\vskip 6pt
Let us now prove the finite covering theorem. If $B-mod=\cup_{i\in I}(B_{i}-mod)$, then $\cup_{i}(B_{i}-mod)\subset \mathcal{L}\subset B-mod$ and $\cup_{i}(B_{i}-mod)=B-mod$, thus $\mathcal{L}=B-mod$. For all $M\in B-mod$, there is an isomorphism $M\backsimeq l(M)$, which proves that the functor $B-mod\rightarrow \prod_{i}B_{i}-mod$ reflects isomorphisms.\vskip 2pt
Reciprocally, if the family of functors $-\otimes_{A}B_{i}$ reflects isomorphisms. For $M\in B-mod$,
\begin{center}
$l(M)\otimes_{A}B_{i}\backsimeq l(M\otimes_{A}B_{i})\backsimeq M\otimes_{A}B_{i}$
\end{center}
Thus $l(M)\backsimeq M$ and $B-mod=\mathcal{L}$.
\vskip 2pt
$\hfill\blacklozenge$

\begin{proposition}\label{R18}
Let $A$ be a commutative monoid and $Y:=spec(B)$ be an affine Zariski open subobject of $X:=Spec(A)$, then the Gabriel filter associated to $B$ is primitive and generated by an ideal of finite type.
\end{proposition}
Proof\\
We write $G_{B}=\cup_{q\in G^{ft}_{B}}G_{q}$ in $GabA$, where $G^{ft}_{B}$ is the set of ideals of finite type in $G_{B}$, and thus $\mathcal{L}_{B}=\cap_{q\in G^{ft}_{B}}\mathcal{L}_{q}$ in $LocA$. In particular $-\otimes_{A}B\backsimeq Colim (l_{q})$ where the functors $l_{q}$ are the left adjoints for the localizations $\mathcal{L}_{q}$. As $B$ is a finitely presentable $A$-algebra, the isomorphism $B\backsimeq colim(l_{q}(A))$ factors through a $l_{q}(A)$. Therefore $l_{q}(A)\in B-mod$, and as all functors $l_{q}$ are isomorphic to identity on $B-mod$, we obtain for all $q'\in G_{B}$ that $l_{q'}\circ l_{q}(A)\backsimeq l_{q}(A)$, hence $B\otimes_{A}l_{q}(A)\backsimeq colim_{q'}(l_{q'}\circ l_{q}(A))\backsimeq l_{q}(A)$. As there is a functorial isomorphism $(B\otimes_{A}-)\circ l_{q}\backsimeq (B\otimes_{A}-)$, $B\backsimeq l_{q}(A)$.\vskip 2pt  
Last but not least, let us prove that $\mathcal{L}_{q}$ is a localization of $B-mod$. The functor $l_{q}$ is left $\mathcal{C}$-adjoint thus commutes with $\mathcal{C}$-colimits. In particular, for all $X\in\mathcal{C}$ and $Y\in A-mod$, $X\otimes Y\in A-mod$ and $l_{q}(X\otimes Y)=X\otimes l_{q}(Y)$. Moreover, if $X,Y$ are two $A$-modules, 
\begin{center}
$X\otimes_{A}Y:=coker(\xymatrix{X\otimes A\otimes Y\ar@<1ex>[r]\ar@<-1ex>[r]&X\otimes Y})$
\end{center}
Thus
\begin{center}
$l_{q}(X\otimes_{A}Y)\backsimeq X\otimes_{A}l_{q}(Y)$
\end{center}
And the functor $l_{q}$ makes a $A$-module $M$ into a  $l_{q}(A)$-module (hence a $B$-module) $l_{q}(M)$. Finally $G_{B}=G_{q}$, i.e. $G_{B}$ is primitive and generated by an ideal of finite type.
\vskip 2pt
$\hfill\blacklozenge$

\subsubsection{Zariski Open subobjects}

\begin{definition}
Let $A$ be a commutative monoid in $\mathcal{C}$ and $F$ be a Zariski open subobject of $Spec(A)$. There exists a family of affine Zariski open subobjects $(Spec(B_{i}))_{i\in I}$ of $Spec(A)$ such that $F$ is the image of $\coprod_{i\in I}Spec(B_{i})\rightarrow Spec(A)$. The Gabriel filter associated to $F$ is $G_{F}:=\cap_{i\in I}G_{B_{i}}$.
\end{definition}
\begin{lemme}
This definition does not depend on the choice of the covering of $F$. 
\end{lemme}
Proof\\
If $(Spec(B_{i}))_{i\in I}$ and $(Spec(C_{j}))_{j\in J}$ are two coverings, define $D_{i,j}:=B_{i}\coprod_{A}C_{j}$. The family $(Spec(D_{i,j}))_{(i,j)\in I\times J}$ is also a covering and each $Spec(B_{i})$ (resp $Spec(C_{j})$) is covered by $(Spec(D_{i,j}))_{j\in J}$ resp ($(Spec(D_{i,j}))_{i\in I}$) and finally $\cap_{i\in I}G_{B_{i}}=\cap_{i,j\in I\times J}G_{D_{i,j}}=\cap_{j\in J}G_{C_{j}}$.
\vskip 2pt
$\hfill\blacklozenge$
\begin{theoreme}\label{lpg}
Let $A$ be a commutative monoid in $\mathcal{C}$ and $F$ be a Zariski open subobject of $Spec(A)$. The Gabriel filter associted to $F$ is locally primitive.
\end{theoreme} 
Proof\\
The scheme $F$ is covered by a family $Spec(B_{i})_{i\in I}$. Let $q_{i}$ be the ideal of finite type generating the Gabriel filter $G_{q_{i}}$ associated to $B_{i}$. As the $q_{i}$ are of finite type $G_{F}=\cap_{i\in I}G_{q_{i}}=\cap_{i\in I}G^{q_{i}}$. Thus $G_{F}=G^{\cup_{i\in I}q_{i}}$.
\vskip 2pt
$\hfill\blacklozenge$ 

\subsection{Zariski Open subobjects in Strong Relative Contexts}
We assume in this section that the relative context $\mathcal{C}$ is strong, i.e. that the functor $(-)_{0}$ reflects isomorphisms. By \ref{prgen}, this implies that, for any object $X$, the set $X_{0}$ is a generating family. 

\subsubsection{Affine Zariski Open subobjects}
\begin{theoreme}\label{basis}
Let $A$ be a commutative monoid in $\mathcal{C}$, any Zariski affine open subobject $Y:=Spec(B)$ of $X:=Spec(A)$ has a finite covering by objects $spec(A_{f})$.
\end{theoreme}
Proof\\The Gabriel filter $G_{B}$ associated to $B$ is primitive and generated by an ideal $q$ of finite type. Therefore, as $\mathcal{C}$ is strong, there exists a generating finite family of elements of $q_{0}$ in $q$. We write $G_{B}=G_{q}=\cap_{f\in F}G_{f}$. As $F$ is finite, \ref{CTF} implies that $(Spec(A_{f}))_{f\in F}$ is a covering of $Spec(B)$.
\vskip 2pt
$\hfill\blacklozenge$

\begin{theoreme}Covering theorem\\
Let $A$ be a commutative monoid and $(Spec(B_{i}))_{i\in I}$ be a family of Zariski open subobjects of $Spec(A)$. Then this family is a covering of the Zariski open subobject $spec(B)$ for the topology defined on $Comm(\mathcal{C})$ (cf \cite{8}) if and only if $B-mod=\cup_{i\in I}B_{i}-mod$ in the poset $LocA$.
\end{theoreme}
Proof\\
Assume that $B-mod=\cup_{i\in I}B_{i}-mod$. We have $G_{q}:=G_{B}=\cap_{i\in I}G_{q_{i}}$ where $G_{q_{i}}:=G_{B_{i}}$. The ideal $q$ is the union of the $q_{i}$, by \ref{prgen} $q$ is isomorphic to the filtered colimit $colim_{J\;fini\subset I}\cup_{j\in J}q_{j}$. As it is of finite type, there exists $J$ such that this isomorphism factors through $q\rightarrow \cup_{j\in J}q_{j}$. The morphism $\cup_{j\in J}q_{j}\rightarrow q$ is then a monomorphism and a split epimorphism hence an isomorphism. By \ref{CTF} the family $(Spec(B_{j}))_{j\in J}$ is a covering of $Spec(B)$ thus so is the family $(Spec(B_{i}))_{i\in I}$.\vskip 2pt 
Reciprocally, if the family $(Spec(B_{i}))_{i\in I}$ is a covering of $Spec(B)$, it has a finite covering sub-family $(Spec(B_{j}))_{j\in J}$ and by \ref{CTF} $B-mod=\cup_{j\in J}B_{j}-mod$. As $\cup_{j\in J}B_{j}-mod\subset \cup_{i\in I}B_{i}-mod\subset B-mod$, we have also $\cup_{i\in I}B_{i}-mod=B-mod$.
\vskip 2pt 
$\hfill\blacklozenge$

\subsubsection{Zariski Open subobjects Classification}

\begin{theoreme}\label{bij}
Let $A$ be a commutative monoid. There is a (contravarient) bijection between the poset of Zariski open subobject of $X:=Spec(A)$, denoted $Ouv(X)$ and the poset of locally primitive Gabriel filters of $A$, denoted $lp(A)$. Moreover a Zariski open subobject of $Spec(A)$ is quasi-compact if and only if its associated Gabriel filter is quasi-compact (and primitive).
\end{theoreme}
Proof\\
By \ref{lpg}, there is a morphism in sets from $Ouv(A)$ to $lp(A)$. We prove first that this morphism is a contravarient morphism of posets. Let $F\subset F'$ be two Zariski open subobjects of $X$, and $(Spec(C_{j}))_{j\in J}$,$(Spec(B_{i}))_{i\in I}$ be respectively coverings of $F'$ and $F$. Let $q$, (resp $q_{j},q_{i}$) be the ideal generating the Gabriel filter associated to $F$ (resp $C_{j},B_{i}$). The scheme $F$ is covered by the family $(Spec(B_{i}\otimes_{A}C_{j}))_{i,j\in I\times J}$. Thus $G_{F}=\cap_{j\in J}\cap_{i\in I}G^{q_{i}.q_{j}}= \cap_{j\in J}G^{q.q_{j}}\supset \cap_{j\in J}G_{q_{j}}=G_{F'}$.\vskip 2pt Let us prove the surjectivity. Let $G^{q}$ be a locally primitive Gabriel filter. As $\mathcal{C}$ is strong, $G^{q}=\cap_{f\in q_{0}}G_{f}$. The Gabriel filter associated to the Zariski open subobject of $X$ defined as the image of $\coprod_{f\in q_{0}}Spec(A_{f})\rightarrow A$ is $G^{q}$.\vskip 2pt
Let us prove the injectivity. Let $F$, $F'$ be two Zariski open subobjects of $X$ with the same associated Gabriel filter. Let $(Spec(A_{f})_{f\in J}$ and $(Spec(A_{f}))_{f\in J'}$ be respectively coverings of $F$ and $F'$. If $q$ and $q'$ are the ideals generated by $J$ and $J'$,  $G^{q}=G_{F}=G_{F'}=G^{q'}$. For all $f\in q_{0}$, there exists $n$ such that $f^{n}\in q'$. As $A_{f^{n}}\backsimeq A_{f}$ any affine Zariski open subobject $A_{f}$ from the covering family of $F$ is a Zariski open subobject of $F'$. Therefore $F$ is a Zariksi open subobject of $F'$. As the roles of $F$ and $F'$ are symmetric, $F\backsimeq F'$.\vskip 2pt
Finally, if $G^{q}$ is a quasi-compact locally primitive Gabriel filter, it is primitive as $G^{q}=G_{q}$. Moreover we can write it as a finite intersection of filters $G_{f}$ and its associated Zariski open subobject is covered by a finite family $Spec(A_{f})$ thus is quasi-compact.
\vskip 2pt
$\hfill\blacklozenge$
\begin{remarque}
The poset $lp(A)$ is a colocale. It is the colocale of the closed subsets of the the Zariski topological space associated to $Spec(A)$.
\end{remarque}
\begin{lemme}\label{cloqc}Let $A$ be in $Comm(\mathcal{C})$. 
\begin{list}{$\triangleright$}{}
\item A locally primitve Gabriel filter of $A$ is irreducible if and only if it has a unique prime representant. 
\item Two irreducible locally primitive Gabriel filter are equals if and only if they have the same subsets of prime ideals.
\item Two locally primitive Gabriel filters are equals if and only if thet have the same subsets of prime ideals.
\end{list}
\end{lemme}
Proof
\begin{list}{-}{}
\item Let $G^{q}$ be an irreducible locally primitive Gabriel filter. The radical $\sqrt{q}$ of $q$ (the intersection of the prime ideals containing $q$) represents the same locally primitive Gabriel filter and is prime. Indeed $q\subset \sqrt{q}$ anq $q$ contains a power of any element of $\sqrt{q}$. If $\sqrt{q}\backsimeq q'.q''$, $q'$ and $q''$ contains $\sqrt{q}$ and as $G^{\sqrt{q}}$ is irreducible, $\sqrt{q}$ contains  a power of any element of one of these two ideals, for exemple $q'$. Then, by definition of the radical of an ideal $\sqrt{q}$ contains $q'$ and thus $\sqrt{q}=q'$. This proves the existence of a prime representant, the unicity is provided by \ref{prime}.
\item Two irreducible locally primitive Gabriel filters which have the same prime ideals are generated by the same prime ideal thus are equal.
\item Let $G^{q},G^{q'}$ be two locally primitive Gabriel filters. We write them as unions of irreducibles ones $G^{q}=\cup_{i\in I}G^{p_{i}}$ and $G^{q'}=\cup_{j\in J}G^{p_{j}}$. For all $i$, there exists $j$ such that $p_{j}\subset p_{i}$ hence $G^{p_{i}}\subset G^{p_{j}}$. Therefore $G^{q}\subset G^{q'}$. As the roles are symmetric, There is equality.
\end{list}
\vskip 2pt 
$\hfill\blacklozenge$

\begin{definition}
The Zariski topological space associated to $X:=spec(A)$, denoted $Lp(A)$ is the topological space whose points are prime ideals of $A$ and in which a closed subset is the set of prime ideals contained in a locally primite Gabriel filter. 
\end{definition}

\begin{theoreme}\label{spac}
Let $A$ be a commutative monoid. The Zariski topological space $Lp(A)$ associated to $A$ is sober and its associated locale is $Ouv(Spec(A))$.
\end{theoreme}
Proof\\By \ref{cloqc} and \ref{bij}.

\subsection{Exemples of Contexts}
\begin{proposition}
Let $A$ be a ring. There is an homeomorphism between the Zariski topological space $Lp(A)$ and the underlaying Zariski topological space of $Spec(A)$.
\end{proposition}
Proof\\By construction.
\vskip 2pt
$\hfill\blacklozenge$
\begin{proposition}
Let $A$ be a differential ring. There is an homeomorphism between the Zariski topological space $Lp(A)$ and the underlaying Zariski topological space of $Spec(A)$ (given in \cite{11}).
\end{proposition} 
Proof\\By construction.
\vskip 2pt
$\hfill\blacklozenge$
\begin{proposition}
Let $A$ be a monoid in set. There is an homeomorphism between the Zariski topological space $Lp(A)$ and the underlaying Zariski topological space of $Spec(A)$ (defined in \cite{10}).
\end{proposition}
Proof:\\The Zariski topological space constructed in \cite{10} is isomorphic to $Lp(A)$.
\vskip 2pt
$\hfill\blacklozenge$

\stepcounter{section}
 
\setcounter{secnumdepth}{-1}
\section{Appendix - Enriched Category Theory}
\small{See \cite{5},\cite{6} and \cite{9} for more details. Let $\mathcal{C}$ be a relative context. We will consider enriched categories over $\mathcal{C}$, refered to as $\mathcal{C}-categories$. Any usual notion will have its enriched equivalent refered to as the $\mathcal{C}-notion$, for example $\mathcal{C}-limit$, $\mathcal{C}-colimit$...  There is an ''underlying category'' functor. To an enriched category, is associated the category whose objects are the same and morphisms between objects procced from the ''underlying set functor'' $Hom_{\mathcal{C}}(1,-)$. 
\begin{definition}
\begin{enumerate}
\item A $\mathcal{C}-category$ $\mathcal{B}$ consists of a collection of objects $ob(\mathcal{B})$, for any two objects $(X,Y)$, of an object $\underline{Hom}_{\mathcal{B}}(X,Y)\in \mathcal{C}$ and of a composition law:
\begin{center}
$M_{X,Y,Z}:\underline{Hom}_{\mathcal{B}}(Y,Z)\otimes\underline{Hom}_{\mathcal{B}}(X,Y)\rightarrow\underline{Hom}_{\mathcal{B}}(X,Z)$
\end{center}
such that the appropriate commutative diagrams commute.
\item A $\mathcal{C}-functor$ $F:\mathcal{B}\rightarrow B'$ consists of, for each object $X$ in $B$, an object $F(X)$ in $B'$, and for each couple of objects $(X,Y)$, a morphism:
\begin{center}
$\xymatrix{\underline{Hom}_{B}(X,Y)\ar[r]^{F_{X,Y}}&\underline{Hom}_{B'}(X,Y)}$
\end{center}
such that the appropriate diagrams commute.
\item An equivalence of $\mathcal{C}-category$ consist of a $\mathcal{C}-functor$ $\mathcal{C}$-fully faithfull such that its underlying functor is essentially surjective.
\end{enumerate}
\end{definition} 
\begin{definition}\label{a2}
\begin{enumerate}
\item Let $\mathcal{B}$ be a $\mathcal{C}-category$, $\mathcal{K}$ a small $\mathcal{C}$-category, and  $F:\mathcal{K}\rightarrow \mathcal{C}$ and $G:\mathcal{K}\rightarrow \mathcal{B}$ two $\mathcal{C}$-functors. Then, if the functor $\mathcal{B}^{op}\rightarrow \mathcal{C}$:
\begin{center}
$B\rightarrow \underline{Hom}_{\mathcal{C}^{\mathcal{K}}}(F,\underline{Hom}_{\mathcal{B}}(B,G-))$
\end{center}
is represented by an object $lim_{F}(G)$, $G$ has a $\mathcal{C}-limit$ indexed by $F$.
\item Let $\mathcal{B}$ be a $\mathcal{C}-category$, $\mathcal{K}$ a small $\mathcal{C}$-category, and  $F:\mathcal{K}^{op}\rightarrow \mathcal{C}$ and $G:\mathcal{K}\rightarrow \mathcal{B}$ two C-functors. Then, if the functor $\mathcal{B}\rightarrow \mathcal{C}$:
\begin{center}
$B\rightarrow \underline{Hom}_{\mathcal{C}^{\mathcal{K}^{op}}}(F,\underline{Hom}_{\mathcal{B}}(G-,B))$
\end{center}
is represented by an object $colim_{F}(G)$, $G$ has a $C-colimit$ indexed by $F$.
\item A $\mathcal{C}-category$ is finite if and only if its collection of objects is a finite set and for each couple of objects $(X,Y)$, the object $\underline{Hom}(X,Y)$ is finitely presented in $\mathcal{C}$.
\item A $\mathcal{C}-limit$ is finite if and only if it is indexed by a functor $F:\mathcal{K}\rightarrow \mathcal{C}_{0}$ from a finite $\mathcal{C}-category$ $\mathcal{K}$ to the generating subcategory of $\mathcal{C}$ consisting of finitely presented objects. 
\end{enumerate}
\end{definition}
For tensored and cotensored (exponential objects exists) $\mathcal{C}$-category , the enriched limit of $G$ along $F$ is in fact the end of $G^{F}$ and the enriched colimit is the coend of $F\otimes G$. Relative contexts are tensored and cotensored self enriched categories.
\begin{proposition}\label{a3} (cf \cite{5}, 3.10)\\
A $\mathcal{C}-Category$ $\mathcal{B}$ is $\mathcal{C}-complete$ if and only if it is Complete and cotensored. Furthermore, $\mathcal{B}$ has all finite $\mathcal{C}-limits$ if and only if it all finite limits and all exponential objects $X^{Y}$, $Y$ finitely presented in $\mathcal{C}$, are representable in $\mathcal{C}$. In particular, a functor $F$ is left $\mathcal{C}-exact$ if and only if it is left exact and commutes with exponential objects $X^{Y}$ (i.e. $F(X^{Y})\backsimeq F(X)^{Y}$), $Y$ finitely presented in $\mathcal{C}$. 
\end{proposition}

\begin{corollaire}\label{a1}
Let $A$ be in $Comm(\mathcal{C})$, then filtered colimits are $\mathcal{C}-exact$ in $A-mod$.
\end{corollaire}
\setcounter{secnumdepth}{2}




\normalsize
\begin{thebibliography}{12}
\bibitem[B]{1} Barr - Exact categories - \textit{Springer, Lecture notes in mathematics 236,} pp1-120, 1971. 
\bibitem[BQ]{2} F.Borceux, C. Quinteiro - A theory of enriched sheaves - \textit{Cahiers de topologie et g\'eom\'etrie diff\'erentielle cat\'egorique} volume XXXVII-2 pages 145-162. 
\bibitem[D]{10} A. Deitmar - Schemes over $F_{1}$ - pr\'e-publication math.NT/0404185
\bibitem[G]{3} P. Gabriel - Des Catégories ab\'eliennes - \textit{Bulletin de la société mathématique de france,} 90, 1962, p 323-448. 
\bibitem[SGA4]{4} M. Artin, A.Grothendieck et J.L. Verdier - Th\'eorie des topos et cohomologie \'etale des sch\'emas. Tome 1: Th\'eorie des topos - S\'eminaire de G\'eom\'etrie Alg\'ebrique du Bois-Marie 1963-1964 (SGA4) - lecture notes in mathematics 269 - \textit{Springer-Verlag, Berlin-New York, 1972}. xix+525pp.
\bibitem[H]{9} M. Hovey - Model Categories - Mathematical Surveys and Monographs, 63 - \textit{American Mathematical Society, Providence, RI,} 1999. xii+209pp.
\bibitem[Kap]{11} I. Kaplanski - An introduction to differential algebra - Hermann, Paris, 1957.
\bibitem[K]{5} G.M. Kelly - Basic concepts of enriched category theory - London Mathematical Society Lecture Note Series, 64 - \textit{Cambridge University Press, Cambridge-New York,} 1982. 285pp. Also available in Reprints in theory and applications of categories, No.10, 2005. 
\bibitem[KD]{6} G.M.Kelly, B.J. Day - Enriched functors categories - 1969 \textit{Reports of the Midwest Category Seminar, III}. pp178-191.  
\bibitem[McL]{7} S. Mac Lane - Categories for the working mathematician - Graduate text in mathematics, 5 - \textit{Springer-Verlag, New York-Berlin , 1971.} ix+262pp. 
\bibitem[TV]{8} B. Toën, M. Vaquie - Under spec(Z) - pr\'e-publication math/O5O9684.

\end{thebibliography}
\end{document}